\documentclass[12pt]{amsart}
\usepackage{amssymb}
\usepackage[all]{xy}

\newcommand{\issuenumber}{9}
\newcommand{\issuemonth}{June}
\newcommand{\issueyear}{2004}

\newtheorem{thm}{Theorem}[section]
\newtheorem{prob}[thm]{Problem}

\newtheorem{issue}{Issue}

\theoremstyle{definition}

\theoremstyle{remark}

\newcommand{\x}{\times}

\newcommand{\Cantor}{{{}^\N\{0,1\}}}

\newcommand{\Iff}{\Leftrightarrow}

\newcommand{\roth}{P_\infty(\N)}

\newcommand{\gpbl}{{\mbox{\textit{\tiny gp}}}}

\renewcommand{\b}{{\mathfrak b}}

\newcommand{\p}{{\mathfrak p}}

\newcommand{\NON}{{\mathsf   {NON}}}
\newcommand{\COF}{{\mathsf   {COF}}}

\newcommand{\dnannouncement}[1]{[\S\ref{#1} below]}

\newcommand{\M}{\mathcal{M}}

\newcommand{\CH}{the Continuum Hypothesis}

\renewcommand{\b}{\mathfrak{b}}
\renewcommand{\t}{\mathfrak{t}}
\newcommand{\h}{\mathfrak{h}}

\renewcommand{\split}{\mathsf{Split}}
\newcommand{\bq}{\begin{quote}}
\newcommand{\eq}{\end{quote}}
\renewcommand{\O}{\mathcal{O}}
\newcommand{\B}{\mathcal{B}}
\newcommand{\BG}{\B_\Gamma}

\newcommand{\sone}{\mathsf{S}_1}    \newcommand{\sfin}{\mathsf{S}_{fin}}
\newcommand{\ufin}{\mathsf{U}_{fin}}

\newcommand{\nin}{\not\in}

\newcommand{\cU}{\mathcal{U}}
\newcommand{\cV}{\mathcal{V}}
\newcommand{\cW}{\mathcal{W}}
\newcommand{\fU}{\mathfrak{U}}
\newcommand{\fV}{\mathfrak{V}}

\newcommand{\naturals}{{\mathbb N}}
\newcommand{\N}{\naturals}
\newcommand{\sm}{\setminus}
\newcommand{\sbst}{\subseteq}
\newcommand{\by}[2]{\par\hfill\emph{#1}, #2}
\newcommand{\nby}[1]{\par\hfill\emph{#1}}
\newcommand{\Tau}{\mathrm{T}}
\newcommand{\CE}{\textsc{CE}}

\newcommand{\be}{\begin{enumerate}}
\newcommand{\ee}{\end{enumerate}}
\newcommand{\bi}{\begin{itemize}}
\newcommand{\ei}{\end{itemize}}
\renewcommand{\i}{\item}

\newcommand{\general}{\small\vfill\par\noindent\hrulefill\par
\noindent\textbf{Previous issues.}
The first issues of this bulletin,
which contain general information (first issue),
basic definitions, research announcements, and open problems (all issues) are available online,
on \arx{math.GN/$x$}, where $x$ is \texttt{0301011}, \texttt{0302062}, \texttt{0303057},
\texttt{0304087}, \texttt{0305367}, \texttt{0312140}, and \texttt{0401155},
\texttt{0403369}, respectively, for issues number $1$ to $8$.\\[0.1cm]
\textbf{Contributions.}
Please submit your contributions (announcements, discussions, and open problems)
by e-mailing us. It is preferred to write them
in \LaTeX{}.
The authors are urged to use as standard notation as possible, or otherwise give
the definitions or a reference to where the notation is explained.
Contributions to this bulletin would not require any transfer of copyright,
and material presented here can be published elsewhere.\\[0.1cm]
\textbf{Subscription.}
To receive this bulletin (free) to your
e-mailbox, e-mail us:\\
{tsaban@math.huji.ac.il}
}

\newcommand{\ArxPaper}[5]{\subsection{#3}{#5}\par\hfill{\arx{#1}}\par\hfill\emph{#4}\par\hfill{#2}}
\newcommand{\nArxPaper}[5]{\subsection{#2}{#4}\par\hfill{\arx{#1}}\par\hfill\emph{#3}}
\newcommand{\AMSPaper}[5]{\subsection{#3}{#5}\par\hfill{\texttt{#1}}\par\hfill\emph{#4}\par\hfill{#2}}
\newcommand{\nAMSPaper}[5]{\subsection{#2}{#4}\par\hfill{\texttt{#1}}\par\hfill\emph{#3}}
\newcommand{\SPMBul}{\textbf{$\mathcal{SPM}$ Bulletin}}

\newcommand{\arx}[1]{\texttt{http://arxiv.org/abs/#1}}
\newcommand{\url}[1]{\bq\texttt{#1}\eq}
\newcommand{\online}[1]{The paper is available online at \url{#1}}

\newcommand{\probmonth}{\emph{Problem of the month}}

\title[$\mathcal{SPM}$ Bulletin \textbf{\issuenumber} (\issuemonth{} \issueyear)]{%
$\mathcal{SPM}$ Bulletin\\[0.5cm]
Issue number \issuenumber: \issuemonth{} \issueyear \CE{}}

\begin{document}
\maketitle

\tableofcontents

\section{Editor's note}

In addition to the interesting research announcements,
we would like to draw your attention to the conference announced
in \dnannouncement{FotFSV}, in which special sessions will be
devoted to infinite games in topology and in set theory, both of interest
to readers of this bulletin. Marion Scheepers is an invited speaker in this
conference, and I also hope to be able to attend. It would be nice to meet
some of you there (that is, in a non-electronic manner).

\subsection*{Email addresses}
Because of the rapidly growing problem of spam emails,
where the addresses are sometimes found by programs surfing in
the internet, we decided to stop, from now on, giving the email
addresses of contributors. If you wish to contact a specific author,
please email me and I will send his email to you personally.

\medskip

Contributions to the next issue are, as always, welcome.

\medskip

\by{Boaz Tsaban}{tsaban@math.huji.ac.il}

\hfill \texttt{http://www.cs.biu.ac.il/\~{}tsaban}

\section{Research announcements}

\subsection{Proceedings of SPM Workshop}
The proceedings of the \emph{Lecce Workshop on Coverings, Selections and Games in Topology}
(June 2002) are going to be published in
\emph{Note di Matematica}, volume \textbf{22}, no.\ 2 (2004).
\nby{Cosimo Guido}

\nAMSPaper{http://www.ams.org/journal-getitem?pii=S0002-9939-04-07351-4}
{A brief remark on van der Waerden spaces}
{Albin L.\ Jones}
{We demonstrate that Martin's axiom for $\sigma$-centered
notions of forcing implies the existence of a van der Waerden space that is
not a Hindman space.  Our proof is an adaptation of the one given by M.\ Kojman
and S.\ Shelah that such a space exists if one assumes the continuum hypothesis
to be true.}

\subsection{Complete ccc Boolean algebras, the order
sequential topology, and a problem of von Neumann}
Let $B$ be a complete $ccc$ Boolean algebra and let
$\tau_s$ be the topology on $B$ induced by the algebraic
convergence of sequences in $B$.
\be
\item Either there
exists a Maharam submeasure on $B$ or every nonempty open set
in $(B,\tau_s)$ is topologically dense.
\item It is consistent that every weakly distributive
complete $ccc$ Boolean algebra carries a strictly positive
Maharam submeasure.
\item The topological space $(B,\tau_s)$ is sequentially
compact if and only if the generic extension by $B$ does not
add independent reals.
\ee
We also give examples of $ccc$ forcings adding a real
but not independent reals.

This paper seems to extend some of the results announced in Section 2.2 of \SPMBul{} 8.
\nby{B.\ Balcar, T.\ Jech, and T.\ Paz\'ak}

\subsection{Cardinal invariants $\p$, $\t$ and $\h$ and real functions}
A partial order on a family of continuous functions from a topological
space $X$ into $[\omega]^{\omega}$ is defined as follows  $$f
\subseteq^{*} g \iff f(x) \subseteq^{*}g(x) \mbox{ for any } x\in X.$$
For these orders variants of cardinals $\mathfrak{p}$, $\mathfrak{t}$ and
$\mathfrak{h}$ are defined and their values are estimated.
\nby{Micha\l{} Machura}

\nArxPaper{math.LO/0404220}
{A comment on $\p<\t$}
{Saharon Shelah}
{We prove that $\p < \t$ if, and only if,
$({}^\omega \omega,<^*)$ has a peculiar cut.\footnote{The
definition of this ``peculiar cut'' appears in the paper.}
We give a self-contained proof (except using Bell theorem).}

\nArxPaper{math.LO/0404421}
{On squares of spaces and $F_\sigma$-sets}
{Arnold W.\ Miller}
{We show that the \CH{} implies there exists
a Lindel\"of space $X$ such that $X^2$ is the union of two metrizable
subspaces but $X$ is not metrizable.  This gives a consistent solution
to a problem of Balogh, Gruenhage, and Tkachuk. The main lemma is
that assuming the \CH{} there exist disjoint sets
of reals $X$ and $Y$ such that $X$ is Borel concentrated on $Y$,
i.e., for any Borel set $B$ if $Y\sbst B$ then $X\sm B$ is countable,
but $X^2\sm\Delta$ is relatively $F_\sigma$ in $X^2\cup Y^2$.}

\nArxPaper{math.LO/0405092}
{Comparing the uniformity invariants of null sets for different measures}
{Saharon Shelah and Juris Stepr\={a}ns}
{It is shown to be consistent with set theory that the uniformity invariant
for Lebesgue measure is strictly greater than the corresponding invariant for
Hausdorff $r$-dimensional measure where $0<r<1$.}

\subsection{Maximal functions and the additivity of various families of null sets}
It is shown to be consistent with set theory that every set of
reals of size $\aleph_1$ is null yet there are $\aleph_1$ planes
in Euclidean $3$-space whose union is not null. Similar results
are obtained for circles in the plane as well as other geometric
objects. The proof relies on results from harmonic analysis about
the boundedness of certain maximal operators and a measure
theoretic pigeonhole principle.
\nby{Juris Steprans}

\ArxPaper{math.GN/0405311}{}
{How many miles to $\beta\omega$? -- Approximating $\beta\omega$ by
metric-dependent compactifications}
{Masaru Kada, Kazuo Tomoyasu, Yasuo Yoshinobu}
{It is known that the Stone-\v{C}ech compactification $\beta{X}$ of a
non-compact metrizable space $X$ is approximated by the collection of Smirnov
compactifications of $X$ for all compatible metrics on $X$. We investigate the
smallest cardinality of a set $D$ of compatible metrics on the countable
discrete space $\omega$ such that, $\beta{\omega}$ is approximated by
Smirnov compactifications for all metrics in $D$, but any finite subset of $D$
does not suffice. We also study the corresponding cardinality for Higson
compactifications.}

\ArxPaper{math.LO/0405473}{}
{The cardinal characteristic for relative gamma-sets}
{Arnold W.\ Miller}
{For $X$ a separable metric space define $\p(X)$ to be the smallest
cardinality of a subset $Z$ of $X$ which is not a relative $\gamma$-set in $X$,
i.e., there exists an $\omega$-cover of $X$ with no $\gamma$-subcover of $Z$. We give
a characterization of $\p(2^\omega)$ and $\p(\omega^\omega)$ in terms of definable
free filters on $\omega$ which is related to the psuedointersection number $\p$.
We show that for every uncountable standard analytic space $X$, either
$\p(X)=\p(2^\omega)$ or $\p(X)=\p(\omega^\omega)$. We show that the following
statements are each relatively consistent with ZFC:
  (a) $\p=\p(\omega^\omega) < \p(2^\omega)$ and
  (b) $\p < \p(\omega^\omega) =\p(2^\omega)$.
}

\nAMSPaper{http://www.ams.org/journal-getitem?pii=S0002-9939-04-07475-1}
{Uncountable intersections of open sets under CPA$_{\mathrm{prism}}$}
{Krzysztof Ciesielski and Janusz Pawlikowski}
{We prove that the Covering Property Axiom CPA$_{\mathrm{prism}}$,
which holds in the iterated perfect set model, implies the following facts.
\be
\i If $G$ is an intersection of $\aleph_1$-many open sets of a Polish space and $G$
has cardinality continuum, then $G$ contains a perfect set.
\i There exists a subset $G$ of the Cantor set which is an intersection of $\aleph_1$-many
open sets but is not a union of $\aleph_1$-many closed sets.
\ee
The example from the second fact refutes a conjecture of Brendle, Larson, and Todorcevic.
}

\AMSPaper{http://www.ams.org/journal-getitem?pii=S0002-9939-04-07422-2}{}
{Covering $\mathbb R^{n+1}$ by graphs of $n$-ary functions and long linear orderings of Turing degrees}
{Uri Abraham and Stefan Geschke}
{A point $(x_0,\dots,x_n)\in X^{n+1}$ is {\em covered} by a
function $f:X^n\to X$ iff there is a permutation $\sigma$ of $n+1$ such that
$x_{\sigma(0)}=f(x_{\sigma(1)},\dots,x_{\sigma(n)})$.
\par
By a theorem of Kuratowski, for every infinite cardinal $\kappa$ exactly
$\kappa$ $n$-ary functions are needed to cover all of $(\kappa^{+n})^{n+1}$.
We show that for arbitrarily large uncountable $\kappa$ it is consistent that
the size of the continuum is $\kappa^{+n}$ and $\mathbb R^{n+1}$ is covered by
$\kappa$ $n$-ary continuous functions.
\par
We study other cardinal invariants of the $\sigma$-ideal on $\mathbb R^{n+1}$
generated by continuous $n$-ary functions and finally relate the question of
how many continuous functions are necessary to cover $\mathbb R^2$ to the
least size of a set of parameters such that the Turing degrees relative to
this set of parameters are linearly ordered.}

\section{Confrences}

\subsection{Foundations of the Formal Sciences V: Infinite Games}\label{FotFSV}
Rheinische Fried\-rich-Wilhelms-Universit\"at Bonn, Mathematisches Institut,
November 26th to 29th, 2004.
\url{http://www.math.uni-bonn.de/people/fotfs/V/}
Infinite Games have been investigated by mathematicians since the
beginning of the twentieth century and have played a central role in
mathematical logic. However, their applications go far beyond mathematics:
they feature prominently in theoretical computer science, philosophical
Gedankenexperiments, as limit cases in economical applications, and in
many other applications. The conference FotFS V wants to bring together
researchers from the various areas that employ infinitary game techniques
to talk about similarities and dissimilarities of the different approaches
and develop cross-cultural bridges.

We invite all researchers from areas applying infinitary game-theoretic
methods (economists, mathematicians, logicians, philosophers, computer
scientists, sociologists) to submit their papers before September 15th,
2004. Topics will include Games in Algebra and Logic, Games in Higher Set
Theory, Games in Set-Theoretic Topology, Infinite Games and Computer
Science, Infinite Games in Philosophy, Infinite Evolutionary Games,
Machine Games, Game Logics, Infinite Games in the Social Sciences.

Invited Speakers:
\bi
\i Samson Abramsky, Oxford UK
\i Alessandro Andretta, Torino
\i Natasha Dobrinen, State College PA
\i Ian Hodkinson, London UK
\i Kevin Kelly, Pittsburgh PA
\i Hamid Sabourian, Cambridge UK
\i Marion Scheepers, Boise ID
\i Brian Skyrms, Irvine CA
\ei
Organizing and Scientific Committee: Stefan Bold (Bonn / Denton TX),
Boudewijn de Bruin (Amsterdam), Peter Koepke (Bonn), Benedikt L\"owe
(Amsterdam / Bonn, Coordinator), Thoralf R\"asch (Potsdam), Johan van
Benthem (Amsterdam / Stanford).

Coordinating e-mail Address: fotfs@math.uni-bonn.de

\nby{Benedikt Loewe}

\section{Problem of the month}

The problem for this issue can be stated with very few definitions.
We will state it this way, and then describe the general framework where it
arises. Recall that $\cU$ is a \emph{large cover} of a set of reals $X$ if
each element of $X$ is covered by infinitely many members of $\cU$.
Let us tentatively say that $X$ has the \emph{splitting property}
if each large open cover of $X$ can be split into two disjoint large
covers of $X$. The \probmonth{} is:

\begin{prob}\label{splitting}
Is it provable that, for all sets of reals $X$ and $Y$ with
the splitting property, $X\cup Y$ has the splitting property?
\end{prob}

The following discussion is based on \cite{splittability}.
Assume that $\fU$ and $\fV$ are collections of covers of a space $X$.
The following property was introduced in \cite{coc1}.
\bi
\i[$\split(\fU,\fV)$:] Every cover $\cU\in\fU$ can be split
into two disjoint subcovers $\cV$ and $\cW$ which contain elements of $\fV$.
\ei
Then the above-mentioned ``splitting property'' is just $\split(\Lambda,\Lambda)$,
where $\Lambda$ denotes the collection of all
large open covers of the space in question.

Properties of this form are useful in the Ramsey theory of thick covers.
The Hurewicz property $\ufin(\Gamma,\Gamma)$
and Rothberger's property $\sone(\O,\O)$
each implies $\split(\Lambda,\allowbreak\Lambda)$, and that
the Sakai property $\sone(\Omega,\Omega)$
implies $\split(\Omega,\Omega)$ \cite{coc1}.
If all finite powers of $X$ have the Hurewicz property (this is equivalent to $\sfin(\Omega^\gpbl,\Omega)$),
then $X$ satisfies $\split(\Omega,\Omega)$ \cite{coc7}.

If we consider this prototype with $\fU,\fV\in\{\Lambda,\Omega,\Tau,\Gamma\}$
we obtain the following $16$ properties.
\newcommand{\aru}{\ar[r]\ar[u]}
$$\xymatrix{
\split(\Lambda,\Lambda)\ar[r]&\split(\Omega,\Lambda)\ar[r]&\split(\Tau,\Lambda)\ar[r]&\split(\Gamma,\Lambda)\\
\split(\Lambda,\Omega)\aru&\split(\Omega,\Omega)\aru&\split(\Tau,\Omega)\aru&\split(\Gamma,\Omega)\ar[u]\\
\split(\Lambda,\Tau)\aru&\split(\Omega,\Tau)\aru&\split(\Tau,\Tau)\aru&\split(\Gamma,\Tau)\ar[u]\\
\split(\Lambda,\Gamma)\aru&\split(\Omega,\Gamma)\aru&\split(\Tau,\Gamma)\aru&\split(\Gamma,\Gamma)\ar[u]
}$$
But all properties in the last column are trivial in the sense that all sets
of reals satisfy them.
On the other hand, all properties but the top one in the first column imply $\binom{\Lambda}{\Omega}$ and are
therefore trivial in the sense that no infinite set of reals satisfies any of them.
Moreover, the properties $\split(\Tau , \Tau)$, $\split(\Tau , \Omega)$, and $\split(\Tau , \Lambda)$ are
equivalent.
It is also easy to see that $\split(\Omega,\Gamma) \Iff \binom{\Omega}{\Gamma}$, therefore
$\split(\Omega,\Gamma)$ implies $\split(\Lambda,\Lambda)$.
In \cite{splittability} it is proved that no implication can be added to the following diagram,
except perhaps the dotted ones.
(If the dotted implication (1) is true, then so are (2) and (3).)
$$
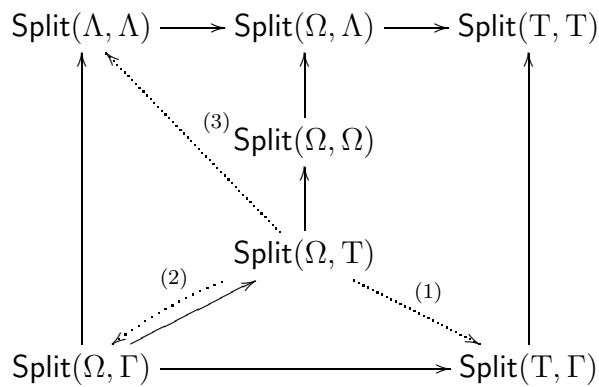
\xymatrix{
\split(\Lambda, \Lambda) \ar[r] & \split(\Omega, \Lambda) \ar[r] & \split(\Tau, \Tau)\\
                           & \split(\Omega, \Omega)\ar[u]\\
& \split(\Omega, \Tau)\ar[u]\ar@{.>}[dr]^{(1)}\ar@/_/@{.>}[dl]_{(2)}\ar@{.>}[uul]_{(3)}\\
\split(\Omega,\Gamma) \ar[uuu]\ar[ur]\ar[rr]     & & \split(\Tau,\Gamma)\ar[uuu]\\
}$$
With regards to the additivity (preservation under taking finite unions) and
$\sigma$-additivity (countable unions), the following is known (\checkmark means that
the property in this position is $\sigma$-additive, and $\x$ means that it is not additive).
$$
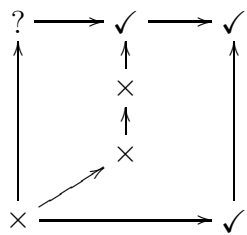
\xymatrix@R=10pt{
\txt{?} \ar[r] & \txt{\checkmark} \ar[r] & \txt{\checkmark}\\
                           & \x\ar[u]\\
& \x\ar[u]\\
\x \ar[uuu]\ar[ur]\ar[rr]     & & \txt{\checkmark}\ar[uuu]\\
}$$

Thus Problem \ref{splitting}, asking whether
$\split(\Lambda,\Lambda)$ is additive, is the only remaining open problem regarding additivity
of these properties.
In Proposition 1.1 of \cite{splittability} it is shown that
for a set of reals $X$ (in fact, for any hereditarily Lindel\"of space $X$),
each large open cover of $X$ contains a \emph{countable} large open cover of $X$.
Consequently, using standard arguments \cite{splittability}, the problem is closely
related to the the following one (where $\roth$ is the space of all infinite
sets of natural numbers, with the topology inherited from $P(\N)$, the latter identified with $\Cantor$).
\begin{prob}
If $\mathsf{R}$ denotes the sets of reals $X$ such that
each continuous image of $X$ in $\roth$ is not reaping,
then is $\mathsf{R}$ additive?
\end{prob}

\by{Boaz Tsaban}{tsaban@math.huji.ac.il}

\section{Problems from earlier issues}
In this section we list the past problems posed in the \SPMBul{},
in the section \probmonth{}.
For definitions, motivation and related results, consult the
corresponding issue.

For conciseness, we make the convention that
all spaces in question are
zero-dimentional, separable metrizble spaces.

\begin{issue}
Is $\binom{\Omega}{\Gamma}=\binom{\Omega}{\Tau}$?
\end{issue}

\begin{issue}
Is $\ufin(\Gamma,\Omega)=\sfin(\Gamma,\Omega)$?
And if not, does $\ufin(\Gamma,\Gamma)$ imply
$\sfin(\Gamma,\Omega)$?
\end{issue}

\begin{issue}
Does there exist (in ZFC) a set satisfying
$\ufin(\O,\O)$ but not $\ufin(\O,\Gamma)$?
\end{issue}
\begin{proof}[Solution]
\textbf{Yes} (Lubomyr Zdomsky).
\end{proof}

\begin{issue}
Does $\sone(\Omega,\Tau)$ imply $\ufin(\Gamma,\Gamma)$?
\end{issue}

\begin{issue}
Is $\p=\p^*$?
\end{issue}

\begin{issue}
Does there exist (in ZFC) an uncountable set satisfying $\sone(\BG,\B)$?
\end{issue}

\begin{issue}
Assume that $X$ has strong measure zero and $|X|<\b$.
Must all finite powers of $X$ have strong measure zero?
\end{issue}
\begin{proof}[Solution]
\textbf{Yes} (Scheepers; Bartoszy\'nski).
\end{proof}

\begin{issue}
Does $X \nin \NON(\M)$ and $Y\nin\mathsf{D}$ imply that
$X\cup Y\nin \COF(\M)$?
\end{issue}

\general

\end{document}